\documentclass{article}
\usepackage{amssymb}
\usepackage{amsfonts}
\usepackage{amsmath}

\newtheorem{theorem}{Theorem}

\newtheorem{lemma}[theorem]{Lemma}

\newtheorem{remark}[theorem]{Remark}

\input{tcilatex}

\begin{document}

\title{Class, degree and integral free forms for the family of Bour's
minimal surfaces}
\author{Erhan G\"{u}ler \\
Bart\i n University, Faculty of Science, \\
Department of Mathematics, Bart\i n, Turkey, \\
Email: ergler@gmail.com}
\maketitle

\begin{abstract}
We consider the family of the Bour's minimal surfaces in Euclidean 3-space,
and compute their classes, degrees and integral free representations.
\end{abstract}

\section{Introduction}

A \textit{minimal surface} in $\mathbb{E}^{3}$ is a regular surface for
which the mean curvature vanishes identically.

Minimal surfaces applicable onto a rotational surface were first determined
by E. Bour \cite{Bour} in 1862. These surfaces have been called $\ \mathfrak{%
B}_{m}$ (following Haag)\ to emphasize the value of $m$. Mathematicians have
dealt with the $\mathfrak{B}_{m}$ in the literature: E. Bour (1862), H.A.
Schwarz (1875), A. Ribaucour (1882), A. Thybaut (1887), A. Demoulin (1897),
L. Bianchi (1899), J. Haag (1906), G. Darboux (1914), E. St\"{u}bler (1914),
J. K. Whittemore (1917), B. Gambier (1921), G. Calugareano (1938).

It was proven by Schwarz \cite{Schwarz} that all real minimal surfaces
applicable to rotational surfaces are given by Whittemore setting%
\begin{equation*}
\mathfrak{F(}s)=Cs^{m-2}
\end{equation*}%
in the Weierstrass representation equations, where $s,C\in \mathbb{C}$, $%
m\in \mathbb{R}$, and $\mathfrak{F(}s)$ is an analytic function. For $C=1,$ $%
m=0$ we obtain the catenoid, $C=i,$ $m=0,$ the right helicoid, $C=1,$ $m=2,$
Enneper's surface (see also \cite{Enneper, Nitsche, Whittemore}). A. Gray 
\cite{Gray} gave the complex forms of the Bour's curve and surface of value $%
m$.

\section{Preliminaries}

Let $\mathbb{E}^{3}$ be three dimensional Euclidean space with natural
metric $\left\langle .\text{ },\text{ }.\right\rangle =dx^{2}+dy^{2}+dz^{2}.$
We will often identify $\overrightarrow{x}$ and $\overrightarrow{x}^{t}$
without further comment.

Let $\mathcal{U}$\ be an open subset of $\mathbb{C}$. A\textbf{\ }\textit{%
minimal\ }(or\textit{\ isotropic})\textit{\ curve}\ is\textit{\ }an analytic
function $\Psi :\mathcal{U}\rightarrow \mathbb{C}^{n}$ such that $\Psi
^{\prime }\left( \zeta \right) \cdot \Psi ^{\prime }\left( \zeta \right) =0,$
where $\zeta \in \mathcal{U}$, and $\Psi ^{\prime }:=\frac{\partial \Psi }{%
\partial \zeta }.$ In addition, if $\Psi ^{\prime }\cdot \overline{\Psi
^{\prime }}=\left\vert \Psi ^{\prime }\right\vert ^{2}\neq 0,$ then $\Psi $
is a \textit{regular minimal curve.} We then have minimal surfaces in the
associated family of a minimal curve, like as given by the following
Weierstrass representation theorem for minimal surfaces.

\begin{theorem}
(K. Weierstrass \cite{Weierstrass}). Let $\mathfrak{F}$ and $\mathcal{G}$ be
two holomorphic functions defined on a simply connected open subset $%
\mathcal{U}$ of $\mathbb{C}$ such that $\mathfrak{F}$ does not vanish on $%
\mathcal{U}$. Then the map%
\begin{equation*}
\mathbf{x}\left( \zeta \right) =\text{Re}\int^{\zeta }\left( 
\begin{array}{c}
\mathfrak{F}\left( 1-\mathcal{G}^{2}\right) \\ 
i\text{ }\mathfrak{F}\left( 1+\mathcal{G}^{2}\right) \\ 
2\mathfrak{F}\mathcal{G}%
\end{array}%
\right) d\zeta
\end{equation*}%
is a minimal, conformal immersion of $\mathcal{U}$ into $\mathbb{%
\mathbb{C}
}^{3},$ and $\mathbf{x}$ is called the Weierstrass patch.
\end{theorem}

\begin{lemma}
Let $\Psi :\mathcal{U}\rightarrow \mathbb{%
\mathbb{C}
}^{3}$ minimal curve and write $\Psi ^{\prime }=\left( \varphi _{1},\varphi
_{2},\varphi _{3}\right) .$ Then%
\begin{equation*}
\mathfrak{F}=\frac{\varphi _{1}-i\varphi _{2}}{2}\text{ \ \ and \ \ }%
\mathcal{G}=\frac{\varphi _{3}}{\varphi _{1}-i\varphi _{2}}
\end{equation*}%
give rise to the Weierstrass representation of $\Psi .$ That is%
\begin{equation*}
\Psi ^{\prime }=\left( \mathfrak{F}\left( 1-\mathcal{G}^{2}\right) ,i%
\mathfrak{F}\left( 1+\mathcal{G}^{2}\right) ,2\mathfrak{F}\mathcal{G}\right)
.
\end{equation*}
\end{lemma}

In section 3, we give the family of Bour's minimal surfaces $\mathfrak{B}%
_{m}.$ We obtain the class and degree of surface $\mathfrak{B}_{3}$ (resp., $%
\mathfrak{B}_{4}$) in section 4 (resp., in section 5). Finally, using the
integral free form of Weierstrass, we find some algebraic functions for $%
\mathfrak{B}_{m}$ $(m\geq 3,$ $m\in 
\mathbb{Z}
)$ in the last section.

\section{The family of Bour's minimal surfaces $\mathfrak{B}_{m}$}

We consider the Bour's curve of value $m.$

\begin{lemma}
The Bour's curve of value $m$%
\begin{equation}
B_{m}\left( \zeta \right) =\left( \frac{\zeta ^{m-1}}{m-1}-\frac{\zeta ^{m+1}%
}{m+1},\text{ }i\left( \frac{\zeta ^{m-1}}{m-1}+\frac{\zeta ^{m+1}}{m+1}%
\right) ,\text{ }2\frac{\zeta ^{m}}{m}\right)  \label{1}
\end{equation}%
is a minimal curve in $\mathbb{%
\mathbb{C}
}^{3}$, where $m\in \mathbb{R}-\left\{ -1,0,1\right\} ,$ $\zeta \in \mathbb{C%
}$, $i=\sqrt{-1}$.
\end{lemma}

We have%
\begin{equation}
B_{m}^{\prime }\cdot B_{m}^{\prime }=0.  \label{2}
\end{equation}%
Bour's surface of value $m$ in $\mathbb{R}^{3}$ is%
\begin{equation}
\mathfrak{B}_{m}\left( \zeta \right) =\text{Re}\int B_{m}^{\prime }\left(
\zeta \right) d\zeta .  \label{3}
\end{equation}

\begin{lemma}
The Weierstrass patch determined by the functions%
\begin{equation*}
\mathfrak{F}\left( \zeta \right) =\zeta ^{m-2}\text{ \ \ and \ \ }\mathcal{G}%
\left( \zeta \right) =\zeta
\end{equation*}%
is a representation of $\mathfrak{B}_{m}.$
\end{lemma}

Therefore, the associated family of minimal surfaces is described by%
\begin{eqnarray*}
\mathfrak{B}\left( r,\theta ;\alpha \right) &=&\text{Re}\int e^{-i\alpha
}B_{m}^{\prime } \\
&=&cos\left( \alpha \right) \text{Re}\int B_{m}^{\prime }+sin\left( \alpha
\right) \text{Im}\int B_{m}^{\prime } \\
&=&cos\left( \alpha \right) \text{ }\mathfrak{B}_{m}\left( r,\theta \right)
+sin\left( \alpha \right) \text{ }\mathfrak{B}_{m}^{\ast }\left( r,\theta
\right) .
\end{eqnarray*}%
When $\alpha =0$ (resp. $\alpha =\pi /2),$ we have the Bour's surface of
value $m$ (resp. the conjugate surface $\mathfrak{B}_{m}^{\ast }$).

The parametric equations of $\mathfrak{B}_{m}$, in polar coordinates $\zeta
=re^{i\theta }$, are%
\begin{equation}
\mathfrak{B}_{m}\left( r,\theta \right) =\left( 
\begin{array}{c}
r^{m-1}\frac{\cos \left[ \left( m-1\right) \theta \right] }{m-1}-r^{m+1}%
\frac{\cos \left[ \left( m+1\right) \theta \right] }{m+1} \\ 
-r^{m-1}\frac{\sin \left[ \left( m-1\right) \theta \right] }{m-1}-r^{m+1}%
\frac{\sin \left[ \left( m+1\right) \theta \right] }{m+1} \\ 
2r^{m}\frac{\cos \left( m\theta \right) }{m}%
\end{array}%
\right) ,  \label{4}
\end{equation}%
with Gauss map%
\begin{equation}
n=\left( \frac{2u}{u^{2}+v^{2}+1},\frac{2v}{u^{2}+v^{2}+1},\frac{%
u^{2}+v^{2}-1}{u^{2}+v^{2}+1}\right) .
\end{equation}

\begin{remark}
$\mathfrak{B}_{m},$ $m\geq 3,$ $m\in 
\mathbb{Z}
,$ has a branch point at $\zeta =0.$ Also, the total curvature of $\mathfrak{%
B}_{3}$ is $-4\pi .$ Note that the catenoid and Enneper's surface are the
only complete regular minimal surfaces in $\mathbb{E}^{3}$ with finite total
curvature $-4\pi $ \cite{Osserman}.
\end{remark}

\begin{remark}
Ribaucour showed that each curve $\mathfrak{B}_{m}\mid _{r=r_{0}}$ lies on
the quadric of revolution%
\begin{equation}
x^{2}+y^{2}+\frac{m^{2}}{m^{2}-1}z^{2}=\left( \frac{r_{0}^{m-1}}{m-1}+\frac{%
r_{0}^{m+1}}{m+1}\right) ^{2}.
\end{equation}
\end{remark}

Next, we will focus on the degree and class of surface $\mathfrak{B}_{m}.$

With $\mathbb{R}^{3}=\{(x,y,z)\mid x,y,z\in \mathbb{R}\},$ the set of roots
of a polynomial $f(x,y,z)=0$ gives an \textit{algebraic surface}. An
algebraic surface is said to be of $degree$ (or $order$) $n$ when $n=\deg
(f).$

The tangent plane on a surface $\mathbf{x}\left( u,v\right) =(x\left(
u,v\right) ,$ $y\left( u,v\right) ,$ $z\left( u,v\right) )$ at a point $%
(u,v) $ is given by%
\begin{equation*}
Xx+Yy+Zz+P=0,
\end{equation*}%
where the Gauss map is $n=(X(u,v),Y(u,v),Z(u,v)),$ $P=P(u,v).$ We have
inhomogeneous tangential coordinates $\overline{u}=X/P,$ $\overline{v}=Y/P,$
and $\overline{w}=Z/P.$ By eliminating $u$ and $v$, we obtain an implicit
equation of $\mathbf{x}\left( u,v\right) $ in tangential coordinates. The
maximum degree of the equation gives $class$ of $\mathbf{x}\left( u,v\right)
.$ See \cite{Nitsche}, for details.

General cases of degree and class of $\mathfrak{B}_{m}$ were studied by
Demoulin \cite{Demoulin}, Haag \cite{Haag}, Ribaucour \cite{Ribaucour} and St%
\"{u}bler \cite{Stubler}. Using the binomial formula we obtain the following
parametric equations of $\mathfrak{B}_{m}\left( u,v\right) :$%
\begin{equation*}
x=\text{Re}\left\{ \frac{1}{m-1}\left[ \sum_{k=0}^{m-1}\tbinom{m-1}{k}%
u^{m-1-k}\left( iv\right) ^{k}\right] -\frac{1}{m+1}\left[ \sum_{k=0}^{m+1}%
\tbinom{m+1}{k}u^{m+1-k}\left( iv\right) ^{k}\right] \right\} ,
\end{equation*}%
\begin{equation*}
y=\text{Re}\left\{ \frac{i}{m-1}\left[ \sum_{k=0}^{m-1}\tbinom{m-1}{k}%
u^{m-1-k}\left( iv\right) ^{k}\right] +\frac{i}{m+1}\left[ \sum_{k=0}^{m+1}%
\tbinom{m+1}{k}u^{m+1-k}\left( iv\right) ^{k}\right] \right\} ,
\end{equation*}%
\begin{equation}
z=\text{Re}\left\{ \frac{2}{m}\left[ \sum_{k=0}^{m}\tbinom{m}{k}%
u^{m-k}\left( iv\right) ^{k}\right] \right\} .
\end{equation}%
It is clear that $\deg \left( x\right) =m+1,$ $\deg \left( y\right) =m+1,$ $%
\deg \left( z\right) =m$ (see also Table 1).

Ribaucour showed that if%
\begin{equation*}
m=\frac{p}{q}\text{ then}\ cl\left( \mathfrak{B}_{m}\right) =2q\left(
p+q\right) ,\ \ m\in 
\mathbb{Z}
\text{ then}\ \deg (\mathfrak{B}_{m})=(m+1)^{2},
\end{equation*}%
\begin{equation*}
m<1\text{ then }cl\left( \mathfrak{B}_{m}\right) =\deg \left( \mathfrak{B}%
_{m}\right) ,\ \ m>1\text{ then }cl\left( \mathfrak{B}_{m}\right) <\deg
\left( \mathfrak{B}_{m}\right) .
\end{equation*}%
Using eliminate methods we calculate the implicit equations, degrees and
classes of the surfaces $\mathfrak{B}_{2},$ $\mathfrak{B}_{3},$ $\mathfrak{B}%
_{4}$. Our findings agree with Ribaucour's, and we give them in Table 1. For
the surface $\mathfrak{B}_{2}$ (i.e., Enneper's surface, see Fig. 1, left
two pictures), it is known that the surface has class $6$, degree $9$. So,
it is also an algebraic minimal surface. For expanded results on $\mathfrak{B%
}_{2}$, see \cite{Nitsche}.

\section{Degree and class of $\mathfrak{B}_{3}$}

The simplest Weierstrass representation $\left( \mathfrak{F},\mathcal{G}%
\right) =\left( \zeta ,\zeta \right) $ gives the \textit{Bour's minimal
surface}\textbf{\ }of value 3. In polar coordinates, the parametric
equations of $\mathfrak{B}_{3}$ (see Fig. 1, right two pictures) are%
\begin{equation}
\mathfrak{B}_{3}\left( r,\theta \right) =\left( 
\begin{array}{c}
\frac{r^{2}}{2}\cos \left( 2\theta \right) -\frac{r^{4}}{4}\cos (4\theta )
\\ 
-\frac{r^{2}}{2}\sin \left( 2\theta \right) -\frac{r^{4}}{4}\sin \left(
4\theta \right) \\ 
\frac{2}{3}r^{3}\cos \left( 3\theta \right)%
\end{array}%
\right) ,
\end{equation}%
where $r\in \lbrack -1,1]$, $\theta \in \lbrack 0,\pi ]$. When $r=1$ on
plane $xy,$ we have \textit{deltoid curve}, which is a 3-cusped hypocycloid
(Steiner's hypocycloid (1856)), also called tricuspoid, is discovered by
Euler in 1745. The parametric form of the surface $\mathfrak{B}_{3},$ in $%
\left( u,v\right) $ coordinates, is%
\begin{equation}
\mathfrak{B}_{3}\left( u,v\right) =\left( 
\begin{array}{c}
-\frac{u^{4}}{4}-\frac{v^{4}}{4}+\frac{3}{2}u^{2}v^{2}+\frac{u^{2}}{2}-\frac{%
v^{2}}{2} \\ 
-u^{3}v+uv^{3}-uv \\ 
\frac{2}{3}u^{3}-2uv^{2}%
\end{array}%
\right) =\left( 
\begin{array}{c}
x(u,v) \\ 
y(u,v) \\ 
z(u,v)%
\end{array}%
\right) ,  \label{6}
\end{equation}%
where $u,v\in \mathbb{R}$. Using the Maple eliminate codes we find the
irreducible implicit equation of surface $\mathfrak{B}_{3}$ as follows:%
\begin{eqnarray*}
&&43046721z^{16}-859963392x^{4}z^{6}-764411904x^{4}y^{2}z^{4} \\
&&-1719926784x^{2}y^{2}z^{6}+509607936x^{2}y^{4}z^{4} \\
&&+69\text{ other lower order terms}=0,
\end{eqnarray*}%
and its degree is $deg(\mathfrak{B}_{3})=16$. Therefore, $\mathfrak{B}_{3}$
is an algebraic minimal surface. To find the class of surface $\mathfrak{B}%
_{3},$ we obtain%
\begin{equation*}
P(u,v)=\frac{(u^{2}+v^{2}+2)(3uv^{2}-u^{3})}{6\left( u^{2}+v^{2}+1\right) },
\end{equation*}%
and the inhomogeneous tangential coordinates%
\begin{eqnarray*}
\overline{u} &=&\frac{12u}{(u^{2}+v^{2}+2)(3uv^{2}-u^{3})}, \\
\overline{v} &=&\frac{12v}{(u^{2}+v^{2}+2)(3uv^{2}-u^{3})}, \\
\overline{w} &=&\frac{6(u^{2}+v^{2}-1)}{(u^{2}+v^{2}+2)(3uv^{2}-u^{3})}.
\end{eqnarray*}%
In tangential coordinates $\overline{u},\overline{v},\overline{w},$ the
irreducible implicit equation of $\mathfrak{B}_{3}$ is%
\begin{eqnarray*}
&&9\overline{u}^{8}+72\overline{u}^{7}+144\overline{u}^{6}+288\overline{u}%
^{5}\overline{w}^{2}+192\overline{u}^{3}\overline{w}^{4}+8\overline{u}^{6}%
\overline{w}^{2} \\
&&-48\overline{u}^{4}\overline{v}^{2}\overline{w}^{2}-576\overline{u}%
\overline{v}^{2}\overline{w}^{4}+81\overline{u}^{2}\overline{v}^{6}+432%
\overline{u}^{4}\overline{v}^{2}-45\overline{u}^{6}\overline{v}^{2} \\
&&-72\overline{u}^{5}\overline{v}^{2}+432\overline{u}^{2}\overline{v}^{4}-360%
\overline{u}^{3}\overline{v}^{4}-216\overline{u}\overline{v}^{6}+27\overline{%
u}^{4}\overline{v}^{4} \\
&&+144\overline{v}^{6}-576\overline{u}^{3}\overline{v}^{2}\overline{w}^{2}+72%
\overline{u}^{2}\overline{v}^{4}\overline{w}^{2}-864\overline{u}\overline{v}%
^{4}\overline{w}^{2}=0.
\end{eqnarray*}%
Therefore, the class of the algebraic minimal surface $\mathfrak{B}_{3}$ is $%
cl(\mathfrak{B}_{3})=8$.

\begin{remark}
Henneberg showed that a plane intersects an algebraic minimal surface in an
algebraic curve \cite{Nitsche}. Using the Gr\"{o}bner eliminate method we
find that the implicit equation of the curve $\mathfrak{B}_{3}\left(
r,0\right) =\mathfrak{\gamma }\left( r\right) =\left( \frac{r^{2}}{2}-\frac{%
r^{4}}{4},0,\frac{2}{3}r^{3}\right) $ (see Fig. $2,$ right two pictures) on
the $xz$-plane is%
\begin{equation*}
1024x^{2}+864xz^{2}-288z^{2}+81z^{4}=0,
\end{equation*}%
and its degree is $deg(\gamma )=4.$ So, we see that the $xz$-plane
intersects the algebraic minimal surface $\mathfrak{B}_{3}$ in an algebraic
curve $\mathfrak{\gamma }\left( r\right) $ (see Fig. $2,$ left two pictures).
\end{remark}

\begin{remark}
The Bour's minimal curve of value $3$ is intersects itself three times along
three straight rays, which meet an angle $2\pi /3$ at the origin in $\mathbb{%
E}^{3}$. The surface $\mathfrak{B}_{3}$ has self-intersections along three
linear rays $u=0,$ $u\pm v\sqrt{3}=0$ at distinct distances from the branch
point $O(0,0,0),$ where $\zeta =u+iv=re^{i\theta }$.
\end{remark}

\section{Degree and class of $\mathfrak{B}_{4}$}

The parametric form of $\mathfrak{B}_{4}$ (see Fig. 3, left two pictures) is%
\begin{equation}
\mathfrak{B}_{4}\left( r,\theta \right) =\left( 
\begin{array}{c}
\frac{r^{3}}{3}\cos \left( 3\theta \right) -\frac{r^{5}}{5}\cos (5\theta )
\\ 
-\frac{r^{3}}{3}\sin \left( 3\theta \right) -\frac{r^{5}}{5}\sin \left(
5\theta \right) \\ 
\frac{1}{2}r^{4}\cos \left( 4\theta \right)%
\end{array}%
\right) ,
\end{equation}%
where $r\in \lbrack -1,1]$, $\theta \in \lbrack 0,\pi ]$. In $\left(
u,v\right) $ coordinates, $\mathfrak{B}_{4}$ has the form as follows%
\begin{equation}
\mathfrak{B}_{4}\left( u,v\right) =\left( 
\begin{array}{c}
\frac{1}{3}u^{3}-uv^{2}-\frac{1}{5}u^{5}+2u^{3}v^{2}-uv^{4} \\ 
-u^{2}v+\frac{1}{3}v^{3}-u^{4}v+2u^{2}v^{3}-\frac{1}{5}v^{5} \\ 
\frac{1}{2}u^{4}-3u^{2}v^{2}+\frac{1}{2}v^{4}%
\end{array}%
\right) ,
\end{equation}%
where $u,v\in \mathbb{R}$. The implicit equation of $\mathfrak{B}_{4}(u,v),$
in cartesian coordinates $x,y,z,$ is as follows 
\begin{eqnarray*}
&&48466299163780426235904z^{25}-147907407116029132800000x^{4}z^{20} \\
&&+887444442696174796800000x^{2}y^{2}z^{20}-147907407116029132800000y^{4}z^{20}
\\
&&-2640558873378816000000000x^{8}z^{15}+233\text{ other lower order terms}=0.
\end{eqnarray*}%
Its degree is $deg(\mathfrak{B}_{4})=25$. Hence, $\mathfrak{B}_{4}$ is an
algebraic minimal surface. To find the class of surface $\mathfrak{B}_{4}$
we obtain%
\begin{equation*}
P(u,v)=\frac{\left( 3u^{2}+3v^{2}+5\right) \left(
v^{4}+6u^{2}v^{2}-u^{4}\right) }{30(u^{2}+v^{2}+1)},
\end{equation*}%
and the inhomogeneous tangential coordinates%
\begin{eqnarray*}
\overline{u} &=&\frac{60u}{\left( 3u^{2}+3v^{2}+5\right) \left(
v^{4}+6u^{2}v^{2}-u^{4}\right) }, \\
\overline{v} &=&\frac{60v}{\left( 3u^{2}+3v^{2}+5\right) \left(
v^{4}+6u^{2}v^{2}-u^{4}\right) }, \\
\overline{w} &=&\frac{30(u^{2}+v^{2}-1)}{\left( 3u^{2}+3v^{2}+5\right)
\left( v^{4}+6u^{2}v^{2}-u^{4}\right) }.
\end{eqnarray*}%
So, the irreducible implicit equation of $\mathfrak{B}_{4},$ in tangential
coordinates $\overline{u},\overline{v},\overline{w},$ is%
\begin{eqnarray*}
&&900\overline{u}^{8}\overline{w}+15\overline{u}^{8}\overline{w}^{2}+15%
\overline{v}^{8}\overline{w}^{2}-180\overline{u}^{2}\overline{w}^{2}-180%
\overline{u}^{2}\overline{v}^{6}\overline{w}^{2}+3600\overline{u}^{2}%
\overline{v}^{6}\overline{w} \\
&&+416\overline{u}^{4}\overline{v}^{6}-3600\overline{u}^{2}\overline{v}%
^{6}-3600\overline{u}^{6}\overline{v}^{2}+8640\overline{u}^{2}\overline{v}%
^{2}\overline{w}^{5}-176\overline{u}^{2}\overline{v}^{8} \\
&&-5400\overline{u}^{4}\overline{v}^{4}+416\overline{u}^{6}\overline{v}%
^{4}-900\overline{v}^{8}\overline{w}-900\overline{v}^{8}+16\overline{v}%
^{10}+16\overline{u}^{10} \\
&&-900\overline{u}^{8}-1440\overline{v}^{4}\overline{w}^{5}-1440\overline{u}%
^{4}\overline{w}^{5}-2400\overline{v}^{6}\overline{w}^{3}+12000\overline{u}%
^{4}\overline{v}^{2}\overline{w}^{3} \\
&&+3600\overline{u}^{6}\overline{v}^{2}\overline{w}-180\overline{u}^{6}%
\overline{v}^{2}\overline{w}^{2}-176\overline{u}^{8}\overline{v}^{2}-2400%
\overline{u}^{6}\overline{w}^{3}-9000\overline{u}^{4}\overline{v}^{4}%
\overline{w} \\
&&+12000\overline{u}^{2}\overline{v}^{4}\overline{w}^{3}+570\overline{u}^{4}%
\overline{v}^{4}\overline{w}^{2}=0.
\end{eqnarray*}%
Hence, the class of the algebraic minimal surface $\mathfrak{B}_{4}$ is $cl(%
\mathfrak{B}_{4})=10$.

We see that the family of $\mathfrak{B}_{m}\left( u,v\right) =(x\left(
u,v\right) ,y\left( u,v\right) ,z\left( u,v\right) )$ are algebraic minimal
surfaces, where $m\in 
\mathbb{Z}
,$ $m\geqslant 2$ (see Table 1).

\section{Integral free form}

Integral free form of the Weierstrass representation is%
\begin{equation}
\left( 
\begin{array}{c}
x \\ 
y \\ 
z%
\end{array}%
\right) =\text{Re}\left( 
\begin{array}{c}
\left( 1-w^{2}\right) \phi ^{\prime \prime }(w)+2w\phi ^{\prime }(w)-2\phi
(w) \\ 
i\left[ \left( 1+w^{2}\right) \phi ^{\prime \prime }(w)-2w\phi ^{\prime
}(w)+2\phi (w)\right] \\ 
2\left[ w\phi ^{\prime \prime }(w)-\phi ^{\prime }(w)\right]%
\end{array}%
\right) \equiv \text{Re}\left( 
\begin{array}{c}
f_{1}\left( w\right) \\ 
f_{2}\left( w\right) \\ 
f_{3}\left( w\right)%
\end{array}%
\right) ,
\end{equation}%
where algebraic function $\phi (w)$ and the functions $f_{i}\left( w\right) $
are connected by the relation%
\begin{equation}
\phi (w)=\frac{1}{4}\left( w^{2}-1\right) f_{1}\left( w\right) -\frac{i}{4}%
\left( w^{2}+1\right) f_{2}\left( w\right) -\frac{1}{2}wf_{3}\left( w\right)
\end{equation}%
for $w\in 
\mathbb{C}
$ \cite{Weierstrass2}$.$ Integral free form is suitable for algebraic
minimal surfaces. For instance, $\phi (w)=\frac{1}{6}w^{3}$ give rise to
Enneper's minimal surface $\mathfrak{B}_{2}$ (see \cite{Nitsche}).

We obtain the function%
\begin{equation}
\phi (w)=\frac{1}{24}w^{4}
\end{equation}%
leads to Bour's minimal surface $\mathfrak{B}_{3}$. We also obtain $\phi _{%
\mathfrak{B}_{4}}(w)=\frac{1}{60}w^{5}$ for $\mathfrak{B}_{4},$ $\phi _{%
\mathfrak{B}_{5}}(w)=\frac{1}{120}w^{6}$ for $\mathfrak{B}_{5},...,$ 
\begin{equation}
\phi _{\mathfrak{B}_{m}}(w)=\frac{1}{(m-1)m\left( m+1\right) }w^{m+1}
\end{equation}
for $\mathfrak{B}_{m},$ where $m\geq 2,m\in 
\mathbb{Z}
.$

\begin{remark}
We find relations between degree of algebraic function $\phi ^{2}(w)$ in the
integral free form and class of surfaces $\mathfrak{B}_{m},$ for integers $%
m\geq 2$. We know $\phi _{\mathfrak{B}_{m}}(w)=\frac{1}{(m-1)m\left(
m+1\right) }w^{m+1}$ for $\mathfrak{B}_{m},$ $m\geq 2,$ $m\in 
\mathbb{Z}
.$ Therefore, we obtain $\deg \left( \phi _{\mathfrak{B}_{2}}^{2}\right)
=6=cl\left( \mathfrak{B}_{2}\right) ,$ $\deg \left( \phi _{\mathfrak{B}%
_{3}}^{2}\right) =8=cl\left( \mathfrak{B}_{3}\right) ,$ $\deg \left( \phi _{%
\mathfrak{B}_{4}}^{2}\right) =10=cl\left( \mathfrak{B}_{4}\right) ,...,$ $%
\deg \left( \phi _{\mathfrak{B}_{m}}^{2}\right) =2m+2=cl\left( \mathfrak{B}%
_{m}\right) .$
\end{remark}

We can see any other parametric eq. and also figure of surface $\mathfrak{B}%
_{m}$ for arbitrary $m\in 
\mathbb{R}
$ using Maple codes. For the figure of $\mathfrak{B}_{5}$ (resp$.$ $%
\mathfrak{B}_{6}$-$\mathfrak{B}_{7},\mathfrak{B}_{8}$-$\mathfrak{B}_{9},%
\mathfrak{B}_{10}$), see Fig. 3, right two pictures (resp. Fig. 4$,$ Fig. 5$%
, $ Fig. 6).

\begin{remark}
We can calculate class of $\mathfrak{B}_{m}$ for integers $m\geq 5,$ but not
calculate degree using Maple codes. Calculation of degree is a time problem
for software programmes.
\end{remark}

\textbf{Acknowledgements. }The author (visiting as a post doctoral
researcher of Katholieke Leuven University, Belgium in 2011-2012 academic
year, and also of Kobe University, Japan at the end of 2013-2014 academic
year) would like to thank the members of the geometry sections, Professor
Franki Dillen (1963-2013), Professor Wayne Rossman, Dr. Ana Irina Nistor and
Masashi Yasumoto for their valuable comments and hospitality.

\begin{center}
\begin{equation*}
\FRAME{itbpF}{4.3777in}{1.3872in}{0in}{}{}{Figure}{\special{language
"Scientific Word";type "GRAPHIC";maintain-aspect-ratio TRUE;display
"USEDEF";valid_file "T";width 4.3777in;height 1.3872in;depth
0in;original-width 9.1774in;original-height 2.885in;cropleft "0";croptop
"1";cropright "1";cropbottom "0";tempfilename
'NHSPQF02.wmf';tempfile-properties "XPR";}}
\end{equation*}

Table 1. Class and degree of $\mathfrak{B}_{m}\left( u,v\right) $ , $m\geq
2, $ $m\in 
\mathbb{Z}
$

\begin{equation*}
\FRAME{itbpF}{4.8516in}{1.2133in}{0in}{}{}{Figure}{\special{language
"Scientific Word";type "GRAPHIC";maintain-aspect-ratio TRUE;display
"USEDEF";valid_file "T";width 4.8516in;height 1.2133in;depth
0in;original-width 11.4164in;original-height 2.8228in;cropleft "0";croptop
"1";cropright "1";cropbottom "0";tempfilename
'NHSPDC00.wmf';tempfile-properties "XPR";}}
\end{equation*}

Figure 1. Left two: Enneper surface $\mathfrak{B}_{2}\left( r,\theta \right) 
$, right two: Bour surface $\mathfrak{B}_{3}\left( r,\theta \right) $

\begin{equation*}
\FRAME{itbpF}{2.7778in}{1.2548in}{0.0138in}{}{}{Figure}{\special{language
"Scientific Word";type "GRAPHIC";maintain-aspect-ratio TRUE;display
"USEDEF";valid_file "T";width 2.7778in;height 1.2548in;depth
0.0138in;original-width 5.6775in;original-height 2.5521in;cropleft
"0";croptop "1";cropright "1";cropbottom "0";tempfilename
'NHSPMR01.wmf';tempfile-properties "XPR";}}
\end{equation*}

Figure 2. Left: Surface $\mathfrak{B}_{3}\left( r,\theta \right) $, right:
its algebraic curve on the xz-plane

\begin{equation*}
\FRAME{itbpF}{4.9675in}{1.1857in}{0in}{}{}{Figure}{\special{language
"Scientific Word";type "GRAPHIC";maintain-aspect-ratio TRUE;display
"USEDEF";valid_file "T";width 4.9675in;height 1.1857in;depth
0in;original-width 11.3956in;original-height 2.6878in;cropleft "0";croptop
"1";cropright "1";cropbottom "0";tempfilename
'NHSPRP03.wmf';tempfile-properties "XPR";}}
\end{equation*}

Figure 3. Left two: Surface $\mathfrak{B}_{4}\left( r,\theta \right) $,
right two: Surface $\mathfrak{B}_{5}\left( r,\theta \right) $

\bigskip

\begin{equation*}
\FRAME{itbpF}{5.0004in}{1.1779in}{0in}{}{}{Figure}{\special{language
"Scientific Word";type "GRAPHIC";maintain-aspect-ratio TRUE;display
"USEDEF";valid_file "T";width 5.0004in;height 1.1779in;depth
0in;original-width 11.2296in;original-height 2.6143in;cropleft "0";croptop
"1";cropright "1";cropbottom "0";tempfilename
'NHSPT804.wmf';tempfile-properties "XPR";}}
\end{equation*}

Figure 4. Left two: Surface $\mathfrak{B}_{6}\left( r,\theta \right) $,
right two: Surface $\mathfrak{B}_{7}\left( r,\theta \right) $%
\begin{equation*}
\FRAME{itbpF}{5.124in}{1.1986in}{0in}{}{}{Figure}{\special{language
"Scientific Word";type "GRAPHIC";maintain-aspect-ratio TRUE;display
"USEDEF";valid_file "T";width 5.124in;height 1.1986in;depth
0in;original-width 11.1353in;original-height 2.5728in;cropleft "0";croptop
"1";cropright "1";cropbottom "0";tempfilename
'NHSPUZ05.wmf';tempfile-properties "XPR";}}
\end{equation*}

Figure 5. Left two: Surface $\mathfrak{B}_{8}\left( r,\theta \right) $,
right two: Surface $\mathfrak{B}_{9}\left( r,\theta \right) $%
\begin{equation*}
\FRAME{itbpF}{2.7216in}{1.1251in}{0in}{}{}{Figure}{\special{language
"Scientific Word";type "GRAPHIC";maintain-aspect-ratio TRUE;display
"USEDEF";valid_file "T";width 2.7216in;height 1.1251in;depth
0in;original-width 5.6559in;original-height 2.3229in;cropleft "0";croptop
"1";cropright "1";cropbottom "0";tempfilename
'NHSPVW06.wmf';tempfile-properties "XPR";}}
\end{equation*}

Figure 6. Left: Surface $\mathfrak{B}_{10}\left( r,\theta \right) ,$ right:
its top view
\end{center}

\end{document}